\documentclass{article}
\usepackage{amsfonts}
\usepackage{amsmath}

\newtheorem{theorem}{Theorem}

\newtheorem{definition}[theorem]{Definition}

\newtheorem{remark}[theorem]{Remark}

\newenvironment{proof}[1][Proof]{\noindent\textbf{#1.} }{\ \rule{0.5em}{0.5em}}
\input{tcilatex}
\begin{document}

\title{Normal and Rectifying Curves in Pseudo-Galilean Space $G_{3}^{1}$ and
Their Characterizations}
\author{Handan \"{O}ztekin \& Alper Osman \"{O}\u{g}renmi\c{s} \and
hbalgetir@firat.edu.tr \& ogrenmisalper@gmail.com \\
F\i rat University, Science Faculty, Mathematics Department TURKEY}
\maketitle

\begin{abstract}
We defined normal and rectifying curves in Pseudo-Galilean Space $G_{3}^{1}.$
Also we obtained some characterizations of this curves in $G_{3}^{1}.$


\end{abstract}

\section{Introduction}

In the Euclidean space $E^{3},$ the notion of rectifying curves was
introduced by B.Y. Chen in [4]. By definition, a regular unit speed space
curve $\alpha (s)$ is called a \linebreak rectifying curve, if its position
vector always lies its rectifying plane $\mathbf{\{t,b\},}$ spanned by the
tangent and the binormal vector field. This subject have been studied by
many researcher. The curves are studied from different way in [4,5,6,7].

A Galilean space may be considered as the limit case of a pseudo-Euclidean
space in which the isotropic cone degenerates to a plane. This limit
\linebreak transition corresponds to the limit transition from the special
theory of \linebreak relativity to classical mechanics. \ [10]

The Pseudo-Galilean \ geometry is one of the real Cayley-Klein geometries
(of projective signature $(0,0,+,-)$. The absolute of the Pseudo-Galilean \
geometry is an ordered triple $\{w,f,I\}$ where \ $w$ is the ideal
(absolute) plane, $f$ is line in $w$ and $I$ is the fixed hyperbolic
involution of points of $f.$ [2]. Differential geometry of the Pseudo -
Galilean space $G_{3}^{1}$ has been largely developed in [1,2,3,8,9].

In the Pseudo-Galilean Space $G_{3}^{1},$ to each regular unit speed curve
\linebreak $r:I\rightarrow G_{3}^{1},$ $I\subset \mathbb{R},$ it is possible
to associate three mutually ortogonalunit vector fields. $\mathbf{t,n}$ and $%
\mathbf{b,}$ called respectively the tangent, the principal normal and the
binormal vector field. The planes spanned by the vector fields $\{\mathbf{t,n%
}\},\{\mathbf{t,b}\}$ and $\{\mathbf{n,b}\}$ are defined as the osculating
plane, the rectifying plane and the normal plane, respectively.

In this paper, we study the normal and rectifying curves in the
Pseudo-Galilean Space $G_{3}^{1}.$ By using similar method as in [4] we show
that there is some characterizations of normal and rectifying curves.

\section{Preliminaries}

Let $r$ be a spatial curve given first by

\begin{equation}
r(t)=(x(t),y(t),z(t)),  \tag{2.1}
\end{equation}%
where $x(t),y(t),z(t)\in C^{3}$ (the set of three-times continuously
differentiable \linebreak functions) and $t$ run through a real interval [2].

\begin{definition}
A curve $r$ given by (2.1) is admissible if
\end{definition}

\begin{equation}
\dot{x}(t)\neq 0.  \tag{2.2}
\end{equation}%
Then the curve $r$ can be given by

\begin{equation}
r(x)=(x,y(x),z(x))  \tag{2.3}
\end{equation}%
and we assume in addition that, in [2]

\begin{equation}
y^{^{\prime \prime }2}(x)-z^{^{\prime \prime }2}(x)\neq 0.  \tag{2.4}
\end{equation}

\begin{definition}
For an admissible curve given by (2.1) the parameter of arc length is
defined by
\end{definition}

\begin{equation}
ds=\left\vert \dot{x}(t)dt\right\vert =\left\vert dx\right\vert .  \tag{2.5}
\end{equation}

For simplicity we assume $dx=ds$ and $x=s$ as the arc length of the curve $r$%
. From now on, we will denote the derivation by $s$ by upper prime $%
^{^{\prime }}$[2].

The vector $\mathbf{t}(s)=r^{^{\prime }}(s)$ is called the tangential unit
vector of an admissible curve $r$ in a point $\mathbf{P}(s)$. Further, we
define the so called osculating plane of $r$ spanned by the vectors $%
r^{^{\prime }}(s)$ and $r^{^{\prime \prime }}(s)$ in the same point. The
absolute point of the osculating plane is

\begin{equation}
H(0:0:y^{^{\prime \prime }}(s):z^{^{\prime \prime }}(s)).  \tag{2.6}
\end{equation}

We have assumed in (2.4) that $H$ is not lightlike. $H$ is a point at
infinity of a line which direction vector is $r^{^{\prime \prime }}(s)$.
Then the unit vector

\begin{equation}
\mathbf{n}(s)=\frac{r^{^{\prime \prime }}(s)}{\sqrt{\left\vert y^{^{\prime
\prime }2}(s)-z^{^{\prime \prime }2}(s)\right\vert }}  \tag{2.7}
\end{equation}%
is called the principal normal vector of the curve $r$ in the point $\mathbf{%
P}$.

Now the vector

\begin{equation}
\mathbf{b}(s)=\frac{(0,\varepsilon z^{^{\prime \prime }}(s),\varepsilon
y^{^{\prime \prime }}(s))}{\sqrt{\left\vert y^{^{\prime \prime
}2}(s)-z^{^{\prime \prime }2}(s)\right\vert }}  \tag{2.8}
\end{equation}%
is orthogonal in pseudo-Galilean sense to the osculating plane and we call
it the binormal vector of the given curve in the point $\mathbf{P}$. Here $%
\varepsilon =+1$ or $-1$ is chosen by the criterion $\det (\mathbf{t},%
\mathbf{n},\mathbf{b})=1$. That means

\begin{equation}
\left\vert y^{^{\prime \prime }2}(s)-z^{^{\prime \prime }2}(s)\right\vert
=\varepsilon (y^{^{\prime \prime }2}(s)-z^{^{\prime \prime }2}(s)). 
\tag{2.9}
\end{equation}%
By the above construction the following can be summarized [2].

\begin{definition}
In each point of an admissible curve in $G_{3}^{1}$ the associated
\linebreak orthonormal (in pseudo-Galilean sense) trihedron $\left\{ \mathbf{%
t}(s),\mathbf{n}(s),\mathbf{b}(s)\right\} $ can be \linebreak defined. This
trihedron is called pseudo-Galilean Frenet trihedron [2].
\end{definition}

If a curve is parametrized by the arc length i.e. given by (2.3), then the
tangent vector is non-isotropic and has the form of

\begin{equation}
\mathbf{t}(s)=r^{^{\prime }}(s)=(1,y^{^{\prime }}(s),z^{^{\prime }}(s)). 
\tag{2.10}
\end{equation}%
Now we have

\begin{equation}
\mathbf{t}^{^{\prime }}(s)=r^{^{\prime \prime }}(s)=(0,y^{^{\prime \prime
}}(s),z^{^{\prime \prime }}(s)).  \tag{2.11}
\end{equation}

According to the clasical analogy we write (2.7) in the form

\begin{equation}
r^{^{\prime \prime }}(s)=\kappa (s)\mathbf{n}(s),  \tag{2.12}
\end{equation}%
and so the curvature of an admissible curve $r$ can be defined as follows

\begin{equation}
\kappa (s)=\sqrt{\left\vert y^{^{\prime \prime }2}(s)-z^{^{\prime \prime
}2}(s)\right\vert }.  \tag{2.13}
\end{equation}

\begin{remark}
In [2] for the pseudo-Galilean Frenet trihedron of an admissible curve $r$
given by (2.3) the following derivative Frenet formulas are true.
\end{remark}

\begin{eqnarray}
\mathbf{t}^{^{\prime }}(s) &=&\kappa (s)\mathbf{n}(s)  \notag \\
\mathbf{n}^{^{\prime }}(s) &=&\tau (s)\mathbf{b}(s)  \TCItag{2.14} \\
\mathbf{b}^{^{\prime }}(s) &=&\tau (s)\mathbf{n}(s)  \notag
\end{eqnarray}%
where $\mathbf{t}(s)$ is a spacelike, $\mathbf{n}(s)$ is a spacelike and $%
\mathbf{b}(s)$ is a timelike vektor, $\kappa (s)$ is the pseudo-Galilean
curvature given by (2.13) and $\tau (s)$ is the pseudo-Galilean torsion of $%
r $ defined by

\begin{equation}
\tau (s)=\frac{y^{^{\prime \prime }}(s)z^{^{\prime \prime \prime
}}(s)-y^{^{\prime \prime \prime }}(s)z^{^{\prime \prime }}(s)}{\kappa ^{2}(s)%
}.  \tag{2.15}
\end{equation}

The formula (2.15) can be written as

\begin{equation}
\tau (s)=\frac{\det (r^{^{\prime }}(s),r^{^{\prime \prime }}(s),r^{^{\prime
\prime \prime }}(s))}{\kappa ^{2}(s)}.  \tag{2.16}
\end{equation}

\section{Normal and Rectifying Curves in Pseudo-Galilean Space $G_{3}^{1}.$}

\begin{definition}
Let $r$ be an admissible curve in 3-dimensional Pseudo-Galilean Space $%
G_{3}^{1}.$ If the position vector of $r$ always lies in its normal plane,
then it is called normal curve in $G_{3}^{1}.$
\end{definition}

By this definition, for a curve in $G_{3}^{1},$ the position vector of $r$
satisfies

\begin{equation}
r(s)=\xi (s)\mathbf{n}(s\mathbf{)+}\eta (s)\mathbf{b}(s),  \tag{3.1}
\end{equation}%
where $\xi (s)$ and $\eta (s)$ are differentiable functions.

\begin{theorem}
Let $r$ be an admissible curve in $G_{3}^{1},$ with $\kappa ,$ $\tau \in 
\mathbb{R}.$ Then $r$ is a normal curve if and only if the principal normal
and binormal components of the position vector are respectively given by
\end{theorem}

\begin{equation}
<r,\mathbf{n}>=(c_{1}+c_{2}s)e^{-\tau s}+(c_{3}+c_{4}s)e^{\tau s}+\frac{%
\kappa }{\tau ^{2}}  \tag{3.2}
\end{equation}%
and

\begin{equation}
<r,\mathbf{b}>=(c_{1}+c_{2}s)e^{-\tau s}-(c_{3}+c_{4}s)e^{\tau s}  \tag{3.3}
\end{equation}%
where $c_{1},c_{2},c_{3},c_{4}\in \mathbb{R}.$

\begin{proof}
Let us assume that $r$ is a normal curve in $G_{3}^{1},$ then from
Definition 1 we have
\end{proof}

\begin{equation}
r(s)=\xi (s)\mathbf{n(}s\mathbf{)+}\eta (s)\mathbf{b}(s).  \tag{3.4}
\end{equation}%
Differentiating this with respect to $s$, we have

\begin{equation}
r^{^{\prime }}(s)=\xi ^{^{\prime }}(s)\mathbf{n(}s\mathbf{)+}\eta ^{^{\prime
}}(s)\mathbf{b}(s)+\xi (s)\mathbf{n}^{\prime }\mathbf{(}s\mathbf{)+}\eta (s)%
\mathbf{b}^{^{\prime }}(s).  \tag{3.5}
\end{equation}%
By using the Frenet equation (2.14), we write

\begin{equation}
\mathbf{t}=\xi ^{^{\prime }}\mathbf{n+}\eta ^{^{\prime }}\mathbf{b}+\xi \tau 
\mathbf{b+}\eta \tau \mathbf{n}.  \tag{3.6}
\end{equation}%
Again differentiating this with respect to $s$ and by using the Frenet
equation (2.14), we get

\begin{equation}
\kappa \mathbf{n}=[(\xi ^{^{\prime }}+\eta \tau )^{^{\prime }}+\tau (\xi
\tau +\eta ^{^{\prime }})]\mathbf{n+[}\tau (\xi ^{^{\prime }}+\eta \tau
)+(\xi \tau +\eta ^{^{\prime }})^{^{\prime }}\mathbf{]b}  \tag{3.7}
\end{equation}%
From equation (3.7), we obtain the differential equation system

\begin{equation}
\begin{array}{c}
\xi ^{^{\prime \prime }}+2\tau \eta ^{^{\prime }}+\tau ^{2}\xi =\kappa \\ 
\eta ^{^{\prime \prime }}+2\tau \xi ^{^{\prime }}+\tau ^{2}\eta =0.%
\end{array}
\tag{3.8}
\end{equation}

By solving this system, we obtain

\begin{equation}
\xi (s)=(c_{1}+c_{2}s)e^{-\tau s}+(c_{3}+c_{4}s)e^{\tau s}+\frac{\kappa }{%
\tau ^{2}},\quad c_{1},c_{2},c_{3},c_{4}\in \mathbb{R}  \tag{3.9}
\end{equation}%
and 
\begin{equation}
\eta (s)=(c_{1}+c_{2}s)e^{-\tau s}-(c_{3}+c_{4}s)e^{\tau s},\quad
c_{1},c_{2},c_{3},c_{4}\in \mathbb{R}  \tag{3.10}
\end{equation}%
which completes the proof.

\begin{definition}
Let $r$ be an admissible curve in 3-dimensional Pseudo-Galilean Space $%
G_{3}^{1}.$ If the position vector of $r$ always lies in its rectifying
plane, then it is called rectifying curve in $G_{3}^{1}.$
\end{definition}

By this definition, for a curve in $G_{3}^{1},$ the position vector of $r$
satisfies

\begin{equation}
r(s)=\lambda (s)\mathbf{t}(s\mathbf{)+}\mu (s)\mathbf{b}(s),  \tag{3.11}
\end{equation}%
where $\lambda (s)$ and $\mu (s)$ are some differentiable functions.

\begin{theorem}
Let $r$ be a rectifying curve in $G_{3}^{1},$ with curvature $\kappa >0,$
\linebreak $<\mathbf{t},\mathbf{t}>=1,$ $<\mathbf{n,n}>=1,$ $<\mathbf{b,b}%
>=\varepsilon ,$ $\varepsilon =\mp 1.$Then the following \linebreak
statements hold:
\end{theorem}

\textit{(i)}The distance function $\rho =\left\Vert r\right\Vert $ satisfies 
\begin{equation*}
\rho ^{2}=\left\vert <r,r>\right\vert =\left\vert
s^{2}+2m_{1}s+m_{1}^{2}+\varepsilon n_{1}^{2}\right\vert 
\end{equation*}
for some $m_{1}\in \mathbb{R},$ $n_{1}\in \mathbb{R}-\{0\}.$

\textit{(ii)} The tangential component of the position vector of $r$ is
given by $<r,\mathbf{t}>\mathbf{=}s+m_{1},$ where $m_{1}\in \mathbb{R}.$

\textit{(iii)} The normal component $r^{N}$ of the position vector of the
curve has a constant length and the distance function $\rho $ is
non-constant.

\textit{(iv) }The torsion $\tau (s)\neq 0$ and binormal component of the
position vector of the curveis constant, i.e. $<r,\mathbf{b}>$ is constant.

\begin{proof}
Let us assume that $r$ is a rectifying curve in $G_{3}^{1}.$ Then from
Definition 3, we can write the position vector of $r$ by%
\begin{equation}
r(s)=\lambda (s)\mathbf{t}(s\mathbf{)+}\mu (s)\mathbf{b}(s),  \tag{3.12}
\end{equation}%
where $\lambda (s)$ and $\mu (s)$ are some differentiable functions of the
ivariant parameters.
\end{proof}

\textit{(i)} Differentiating the equation (3.12) with respect to $s$ and
considering the Frenet equations (2.14), we get 
\begin{equation}
\begin{array}{c}
\lambda ^{^{\prime }}(s)=1 \\ 
\lambda (s)\kappa (s)+\mu (s)\tau (s)=0 \\ 
\mu ^{^{\prime }}(s)=0.%
\end{array}
\tag{3.13}
\end{equation}%
Thus, we obtain

\begin{equation}
\begin{array}{c}
\lambda (s)=s+m_{1},\quad m_{1}\in \mathbb{R} \\ 
\mu (s)=n_{1,}\quad n_{1}\in \mathbb{R} \\ 
\mu (s)\tau (s)=-\lambda (s)\kappa (s)\neq 0,%
\end{array}
\tag{3.14}
\end{equation}%
and hence $\mu (s)=n\neq 0,$ $\tau (s)\neq 0.$From the equation (3.12), we
easily find that

\begin{equation}
\rho ^{2}=\left\vert <r,r>\right\vert =\left\vert
s^{2}+2m_{1}s+m_{1}^{2}+\varepsilon n_{1}^{2}\right\vert ,\quad \varepsilon
=\mp 1  \tag{3.15}
\end{equation}

\textit{(ii)} If we consider equation (3.12), we get

\begin{equation}
<r,\mathbf{t}>\mathbf{=}\lambda (s)  \tag{3.16}
\end{equation}%
which means that the tangential component of the position vector of $r$ is
given by

\begin{equation}
<r,\mathbf{t}>\mathbf{=}s+m_{1},\quad m_{1}\in \mathbb{R}.  \tag{3.17}
\end{equation}

(\textit{iii)} From the equation (3.12), it follows that the normal
component $r^{N}$ of the position vector $r$ is given by 
\begin{equation}
r^{N}=\mu \mathbf{b.}  \tag{3.18}
\end{equation}%
Therefore,%
\begin{equation}
\left\Vert r^{N}\right\Vert =\left\vert \mu \right\vert =\left\vert
n_{1}\right\vert \neq 0.  \tag{3.19}
\end{equation}%
Thus we proved statement \textit{(iii). }

\textit{(iv) }If we consider equation (3.12), we easily get 
\begin{equation}
<r,\mathbf{b}>=\varepsilon \mu =const.,\quad \varepsilon =\mp 1  \tag{3.20}
\end{equation}%
and since $\tau (s)\neq 0$, the statement \textit{(iv) }is proved.

Conversely, suppose that statement \textit{(i)} or statement \textit{(ii)}
holds. Then we have 
\begin{equation}
<r,\mathbf{t}>\mathbf{=}s+m_{1},\quad m_{1}\in \mathbb{R}.  \tag{3.21}
\end{equation}%
Differentiating equation (3.21) with respect to $s$, we obtain%
\begin{equation}
\kappa <r,\mathbf{n}>=0.  \tag{3.22}
\end{equation}%
Since $\kappa >0,$ it follows that%
\begin{equation}
<r,\mathbf{n}>=0  \tag{3.23}
\end{equation}%
which means that $r$ is a rectifying curve.

Next, suppose that statement \textit{(iii) }holds. Let us can write 
\begin{equation}
r(s)=l(s)\mathbf{t}(s\mathbf{)+}r^{N},\quad l(s)\in \mathbb{R}.  \tag{3.24}
\end{equation}%
Then we easily obtain that

\begin{equation}
<r^{N},r^{N}>=C=const.=<r,r>-<r,\mathbf{t}>^{2}.  \tag{3.25}
\end{equation}%
If we differentiate equation (3.25) with respect to $s$, we get%
\begin{equation}
<r,\mathbf{t}>=<r,\mathbf{t}>[1+\kappa <r,\mathbf{n}>].  \tag{3.26}
\end{equation}%
Since $\rho \neq const.,$ we have 
\begin{equation}
<r,\mathbf{t}>\neq 0.  \tag{3.27}
\end{equation}%
Moreover, since $\kappa >0$ and from (3.26) we obtain%
\begin{equation}
<r,\mathbf{n}>=0,  \tag{3.28}
\end{equation}%
that is $r$ is rectifying curve.

Finally, if the statement \textit{(iv)} holds, then from the Frenet
equations (2.14), we get%
\begin{equation}
<r,\mathbf{n}>=0,  \tag{3.29}
\end{equation}%
which means that $r$ is rectifying curve.

\begin{theorem}
Let $r$ be a curve in $G_{3}^{1}.$ Then the curve $r$ is a rectifying curve
if and only if there holds
\end{theorem}

\begin{equation}
\frac{\tau (s)}{\kappa (s)}=as+b  \tag{3.30}
\end{equation}%
where $a\in \mathbb{R}-\{0\},$ $b\in \mathbb{R}.$

Proof. Let us first suppose that the curve $r(s)$ is rectifying. From the
equations (3.13) and (3.14) we easily find that%
\begin{equation}
\frac{\tau (s)}{\kappa (s)}=as+b  \tag{3.31}
\end{equation}%
where $a\in \mathbb{R}-\{0\},$ $b\in \mathbb{R}.$

Conversely, let us suppose that $\frac{\tau (s)}{\kappa (s)}=as+b,$ $a\in 
\mathbb{R}-\{0\},$ $b\in \mathbb{R}.$ Then we may choose

\begin{equation}
\begin{array}{c}
a=\frac{1}{n_{1}} \\ 
b=\frac{m_{1}}{n_{1}}%
\end{array}
\tag{3.32}
\end{equation}%
where $n_{1}\in \mathbb{R}-\{0\},$ $m_{1}\in \mathbb{R}.$

Thus we have

\begin{equation}
\frac{\tau (s)}{\kappa (s)}=\frac{s+m_{1}}{n_{1}}.  \tag{3.33}
\end{equation}

If we consider the Frenet equations (2.14), we easily find that

\begin{equation}
\frac{d}{ds}[r(s)-(s+m_{1})\mathbf{t}(s)-n_{1}\mathbf{b}(s)]=0  \tag{3.34}
\end{equation}%
which means that $r$ is a rectifying curve.

\textbf{REFERENCES}

[1] Divjak, B., Geometrija pseudogalilejevih prostora, Ph.D. thesis,
University of Zagreb, 1997.

[2] Divjak, B., Curves in Pseudo-Galilean Geometry, \textit{Annales Univ.
Sci. Budapest,} 41 (1998) 117-128,

[3] Divjak, B. and Sipus, Z.M., Special curves on ruled surfaces in Galilean
and pseudo-Galilean spaces, Acta Math. Hungar.,98(3) (2003) 203-215.

[4] Chen, B.Y., When does the position vector of a space curve always lie in
its rectifying plane?, Amer. Math. Monthly 110 (2003) 147-152.

[5] Chen, B.Y., Dillen, F., Rectifying curves as centrodes and extremal
curves, Bull. Inst. Math. Academia Sinica, 33(2) (2005) 77-90.

[6] Ilarslan, K., Ne\v{s}ovi\'{c}, E., Petrovi\'{c}-Torga\v{s}ev, M., Some
characterizations of rectifying curves in Minkowski 3-space, Novi Sad J.
Math. 33(2) (2003), 23-32.

[7] Ilarslan, K., Ne\v{s}ovi\'{c}, E., On Rectifying Curves as Centrodes and
Extremal Curves in the Minkowski 3-Space, Novi Sad J. Math. 37(1) (2007),
53-64.

[8] \"{O}\u{g}renmi\c{s}, A.O., Ruled Surfaces in the Pseudo - Galilean
Space, Ph.D. Thesis, University of F\i rat, 2007.

[9] \"{O}\u{g}renmi\c{s}, A.O. and Erg\"{u}t, M., On the Explicit
Characterization of Admissible Curve in 3-Dimensional Pseudo - Galilean
Space, J. Adv. Math. Studies, Vol.2, No.1 (2009) 63-72.

[10] Yaglom, I. M., A Simple Non-Euclidean Geometry and Its Physical Basis,
Springer-Verlag, New York Inc. 1979

\end{document}